\documentclass[12pt]{article}
\usepackage{amsfonts}

\usepackage{latexsym}
\usepackage{amssymb}
\usepackage{amsmath}
\usepackage{amsthm}
\usepackage[mathscr]{eucal}
\usepackage{graphicx}
\usepackage{hyperref}
\usepackage{caption}
\usepackage{subcaption}

\newtheorem{Lemma}{Lemma}

\newtheorem{Theorem}[Lemma]{Theorem}

\newtheorem{Definition}{Definition}

\renewcommand{\qed}{\hfill{\ \ \rule{2mm}{2mm}} \vspace{0.2in}}

\begin{document}

\title{Bridged Hamiltonian Cycles in Sub-critical Random Geometric Graphs}
\author{ \textbf{Ghurumuruhan Ganesan}
\thanks{E-Mail: \texttt{gganesan82@gmail.com} } \\
\ \\
Institute of Mathematical Sciences, HBNI, Chennai}
\date{}
\maketitle

\begin{abstract}
In this paper, we consider a random geometric graph (RGG)~\(G\) on~\(n\) nodes with adjacency distance~\(r_n\) just below
the Hamiltonicity threshold and construct Hamiltonian cycles using additional edges called bridges.
The bridges by definition do not belong to~\(G\) and we are interested in estimating the number of bridges and the maximum bridge length, needed for constructing a Hamiltonian cycle. In our main result, we show that with high probability, i.e. with probability converging to one as~\(n \rightarrow \infty,\) we can obtain a Hamiltonian cycle with maximum bridge length a constant multiple of~\(r_n\) and containing an arbitrarily small fraction of edges as bridges. We use a combination of backbone construction and iterative cycle merging to obtain the desired Hamiltonian cycle.


\vspace{0.1in} \noindent \textbf{Key words:} Random geometric graphs, Hamiltonian cycles with bridges.

\vspace{0.1in} \noindent \textbf{AMS 2000 Subject Classification:} Primary:
60J10, 60K35; Secondary: 60C05, 62E10, 90B15, 91D30.
\end{abstract}

\bigskip

\setcounter{equation}{0}
\renewcommand\theequation{\thesection.\arabic{equation}}
\section{Introduction}\label{intro}
Hamiltonian cycles in Random Geometric Graphs (RGGs) are extremely important from both theoretical and application perspectives. Penrose (1997) obtained sharp bounds on the threshold adjacency distance for the RGG to become Hamiltonian with high probability and later D\'iaz et al (2007) generalized this to metrics other than Euclidean distance along with providing an algorithm for finding the Hamiltonian cycle. Nearly simultaneously, M\"uller et al. (2011), Balogh et al (2011) investigated the problem of coincidence of~\(2-\)connectedness and Hamiltonicity of RGGs and study how the graph becomes Hamiltonian just as it also becomes~\(2-\)connected. More recently Bal et al. (2017) have explored the existence of rainbow Hamiltonian cycles in edge coloured RGGs.

In this paper,  we study construction of Hamiltonian cycles in an RGG~\(G\) with adjacency distance~\(r_n\) just \emph{below} the Hamiltonicity threshold, using a small number of extra edges called bridges that do not belong to~\(G.\) We use dense component construction involving discretization of the unit square, to obtain a number of small cycles in~\(G\) and then ``stitch" the cycles together using bridges to obtain the desired bridged cycle. We also obtain bounds on the adjacency distance that ensures that the bridge fraction in the resulting cycle is arbitrarily small.

The paper is organized as follows. In Section~\ref{sec_rob_rgg}, we describe our main result Theorem~\ref{thm_rob_rgg} regarding the maximum bridge length and bridge fraction of bridged Hamiltonian cycles in random geometric graphs whose adjacency distance is just below the Hamiltonian threshold. Next, in Section~\ref{sec_prelim}, we collect the preliminary results used in the proof of Theorem~\ref{thm_rob_rgg} and finally, in Section~\ref{sec_proof}, we prove Theorem~\ref{thm_rob_rgg}. 


\setcounter{equation}{0}
\renewcommand\theequation{\thesection.\arabic{equation}}
\section{Bridged Hamiltonian Cycles in RGGs}\label{sec_rob_rgg}
Consider $n$ nodes \(X_1,X_2,\ldots,X_n,\) independently distributed in the unit square $S = \left[-\frac{1}{2},\frac{1}{2}\right]^2$ each according to a certain density \(f\) satisfying
\begin{equation}\label{f_eq}
0 < \epsilon_1 \leq f(x) \leq \epsilon_2 < \infty.
\end{equation}
We define the overall process on the probability space \((\Omega_X, {\cal F}_X, \mathbb{P})\) and let~\(K_n = K_n(X_1,\ldots,X_n)\)  be the complete graph with vertex set~\(\{X_1,\ldots,X_n\}.\) Let~\(G(r_n)\) be the graph formed by the set of all edges of~\(K_n,\) each of whose length is strictly less than~\(r_n.\) We define~\(G(r_n)\) to be the random geometric graph (RGG) formed by the nodes~\(\{X_i\}_{1 \leq i \leq n}\) with adjacency distance~\(r_n.\)



A cycle in~\(K_n\) is a sequence of distinct nodes~\({\cal D} = (X_{i_1},\ldots,X_{i_{t}})\) such that~\(X_{i_j}\) is adjacent to~\(X_{i_{j+1}}\) for~\(1 \leq j \leq t-1\) and~\(X_{i_t}\) is adjacent to~\(X_{i_1}.\) The length of~\({\cal D}\) is the number of edges in~\({\cal D}.\) The cycle~\({\cal D}\) is said to be \emph{Hamiltonian} if~\(t = n;\) i.e. the cycle~\({\cal D}\) contains all the~\(n\) nodes~\(\{X_i\}_{1 \leq i \leq n}.\) An edge~\(e \in K_n\) of length at least~\(r_n\) is said to be a \emph{bridge} with respect to~\(G(r_n).\)  Throughout, we consider only bridges with respect to~\(G(r_n)\) and so we suppress the phrase ``with respect to~\(G(r_n)\)". A cycle~\({\cal D} \subset K_n\) is said to be a \emph{bridged cycle} if~\({\cal D}\) contains at least one bridge.
\begin{Definition}\label{def_one} Let~\({\cal D} \subset K_n\) be a Hamiltonian cycle. For~\(w >0\) and~\(0 < \gamma < 1\) we say that~\({\cal D}\) is a~\((w,\gamma)-\)bridged Hamiltonian cycle if the following two properties hold:\\
\((a)\) The maximum length of an edge in~\({\cal D}\) is less than~\(w.\)\\
\((b)\) If~\(n_{br}\) denotes the number of bridges in~\({\cal D},\) then the ratio~\(\frac{n_{br}}{n} \leq \gamma.\)
\end{Definition}
In other words, the fraction of edges that are bridges in~\({\cal D}\) is at most~\(\gamma.\) From the above definition, we see that if~\({\cal D}\) does not contain any bridges, then~\({\cal D}\) is a~\((r_n,0)-\)bridged Hamiltonian cycle.

Given~\(w > 0\) and~\(0 \leq \gamma \leq 1,\) let~\(E(w,\gamma)\) be the event that there exists a~\((w,\gamma)-\)bridged Hamiltonian cycle. We are interested in estimating the probability of occurrence of~\(E(w,\gamma)\) for various values of~\(w\) and~\(\gamma.\) For example, using the fact that any edge of the complete graph~\(K_n\) is at most~\(\sqrt{2},\) we get that~\(E(\sqrt{2},1)\) occurs with probability one. On the other hand, if the nodes are uniformly distributed and the adjacency distance~\(r_n\) is larger than the Hamiltonicity threshold, i.e. if~\(\pi nr_n^2 = \log{n} + \delta \log{\log{n}}\) for some constant~\(\delta > 4,\) then we know (Penrose (1997)) that~\(E(r_n,0)\) occurs with high probability, i.e. with probability converging to one as~\(n \rightarrow \infty.\)

For values of~\(r_n\) just below the Hamiltonicity threshold, we have the following result.
\begin{Theorem}\label{thm_rob_rgg} Let~\(\epsilon_1\) be as in~(\ref{f_eq}) and for a constant~\(\alpha >0,\) define
\begin{equation}\label{rn_hom}
nr_n^2 = \frac{1}{4\epsilon_1}\left(\log{n} + \alpha\log{\log{n}} + \omega_n\right)
\end{equation}
where~\(\omega_n \rightarrow \infty\)  as~\(n \rightarrow \infty.\) For every integer~\(L \geq 9\) and every~\(\alpha \geq 8L-1,\) there exists a constant~\(C >0\) such that
\begin{equation}\label{ew_est}
\mathbb{P}\left(E(2r_n,\gamma)\right) \geq 1-\frac{1}{n^{9}} - \frac{Ce^{-\frac{3\omega_n}{8}}}{(\log{n})^{\alpha-8L+1}},
\end{equation}
where~\(\gamma := \frac{16}{L-8}.\)
\end{Theorem}
Choosing~\(L\) large, the fraction of bridges~\(\gamma\) can be made arbitrarily close to~\(0.\)  If the node distribution is uniform, then~\(\epsilon_1\) as defined in~(\ref{f_eq}) equals one and so for any constant~\(0 < \epsilon  <1\) and all~\(n\) large, we have from~(\ref{rn_hom}) that~\[r_n < \sqrt{\frac{(1-\epsilon)\log{n}}{\pi n}}.\] This implies that the corresponding RGG with adjacency distance~\(r_n\) is in fact below the connectivity threshold and is therefore disconnected with high probability (Penrose (2003)). Theorem~\ref{thm_rob_rgg} says that with high probability, using at most~\(\gamma n\) bridges each of length at most~\(2r_n,\) we can still ``patch" together a Hamiltonian cycle. 

In our proof of Theorem~\ref{thm_rob_rgg} below, we first construct a collection of small cycles formed by edges of~\(G(r_n)\) and then join these cycles together using bridges of length at most~\(2r_n,\) to obtain the desired bridged Hamiltonian cycle.






\setcounter{equation}{0}
\renewcommand\theequation{\thesection.\arabic{equation}}
\section{Preliminaries}\label{sec_prelim}
\begin{figure}[tbp]
\centering
\includegraphics[width=2.8in, trim= 30 290 100 100, clip=true]{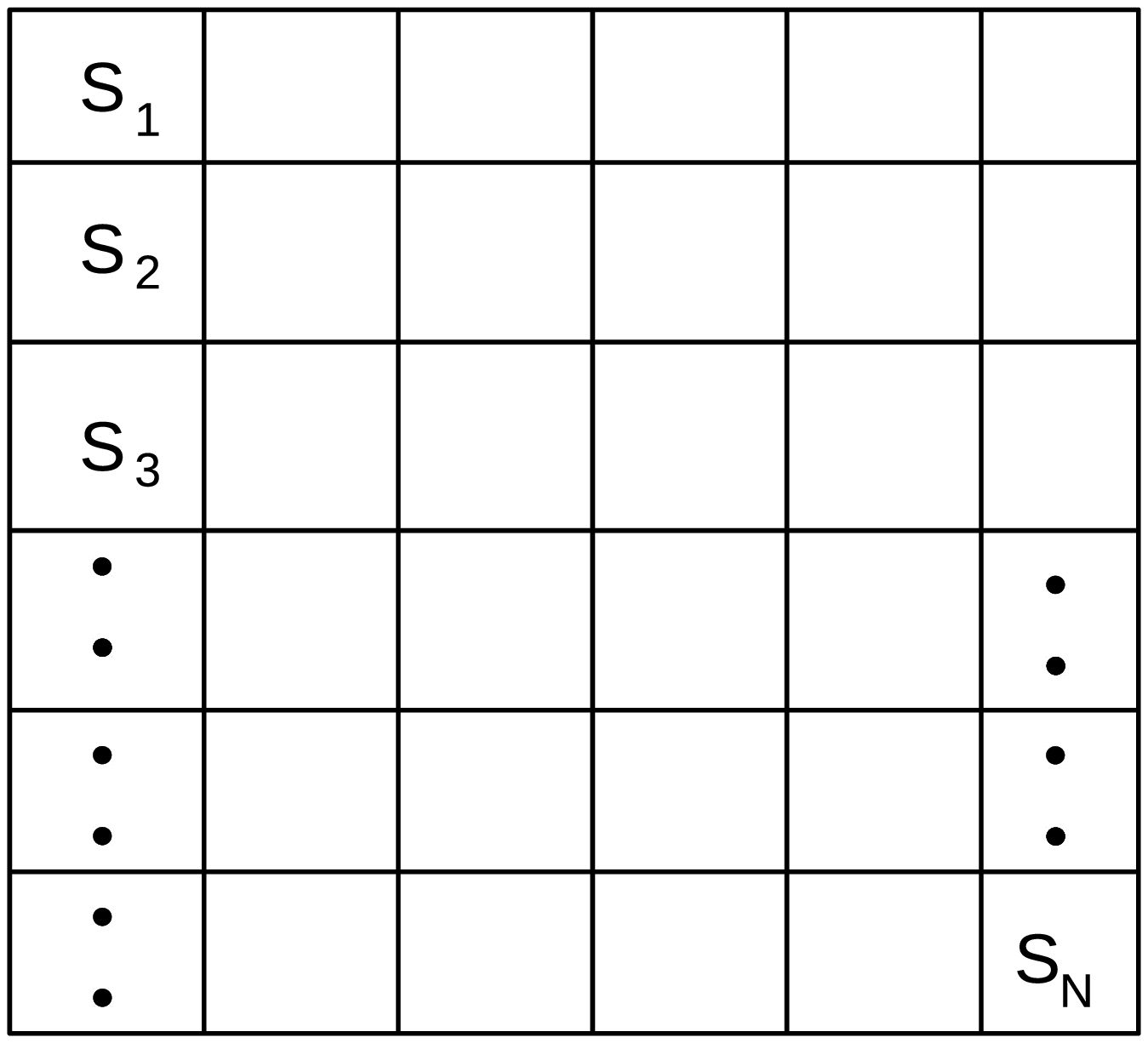}
\caption{Tiling the unit square into disjoint~\(t_n \times t_n\) squares~\(\{S_i\}_{1 \leq i \leq N}.\) }
\label{fig_tile}
\end{figure}

In this section, we collect a couple of preliminary results used in the proof of Theorem~\ref{thm_rob_rgg}. Throughout we use the following discretization procedure.   Divide the unit square~\(S\) into disjoint squares~\(\{S_j\}\) of side length~\(t_n\) as shown in Figure~\ref{fig_tile}, where
\begin{equation}\label{tn_def}
8nt^2_n = \frac{1}{\epsilon_1}\left(\log{n} + \alpha\log{\log{n}} + \omega_n\right) - \gamma_n
\end{equation}
and~\(\gamma_n \in (0,1) \) is such that~\(\frac{1}{t_n}\) is an integer for all~\(n\) large. Such a~\(\gamma_n\) always exists and for completeness we provide a small justification in the Appendix.

By choice~\(t_n\) is slightly less than~\(\frac{r_n}{\sqrt{2}}\) and so any two nodes within a square~\(S_i\) are adjacent in the graph~\(G(r_n).\) Moreover if squares~\(S_{i}\) and~\(S_{j}\) are adjacent (i.e., share a corner), then every node of~\(\{X_k\}\) in~\(S_{i}\) is joined to every node in~\(S_{j}\) by an edge of length less than~\(2r_n.\) Say that~\(S_{j}\) is \emph{dense} if it contains at least~\(L\) nodes of~\(\{X_k\}\) and \emph{sparse} otherwise. Letting~\(E(j)\) be the event that~\(S_j\) is sparse, the following Lemma estimates the joint event that multiple squares are sparse. Throughout constants do not depend on~\(n.\)
\begin{Lemma}\label{good_lem} Let~\(q \geq 1\) be any integer constant and let~\(S_{j_1},\ldots,S_{j_q}\) be any set of~\(q\) squares in~\(\{S_i\}_{i \geq 1}.\) We have that
\begin{equation}\label{sj_sparse}
\mathbb{P}\left(\bigcap_{i=1}^{q}E(j_i)\right) \leq \frac{C}{n^{q/8} \cdot (\log{n})^{q\alpha_1}}\exp\left(-\frac{q\omega_n}{8}\right)
\end{equation}
for all~\(n\) large, where~\(\alpha_1  := \frac{\alpha}{8}-L\) and~\(C  = C(q) > 0\) is a positive constant.
\end{Lemma}
\emph{Proof of Lemma~\ref{good_lem}}: If~\(\bigcap_{i=1}^{q}E(j_i)\) occurs, then the total number of nodes in~\(\bigcup_{i=1}^{q} S_{j_i}\)
is at most~\(Lq.\) The total area covered by~\(\bigcup_{i=1}^{q} S_{j_i}\) is~\(qt_n^2\) and so the probability that node~\(X_i\)
belongs to one of these~\(q\) squares is~\[\int_{\bigcup_{i=1}^{q} S_{j_i}} f(x)dx \in [\epsilon_1 qt_n^2, \epsilon_2 qt_n^2],\] using the bounds for~\(f(.)\) in~(\ref{f_eq}). Thus
\begin{eqnarray}
\mathbb{P}\left(\bigcap_{i=1}^{q}E(j_i)\right) &\leq& \sum_{k=0}^{Lq} {n \choose k} (\epsilon_2 qt_n^2)^{k} (1-\epsilon_1qt_n^2)^{n-k} \nonumber\\
&\leq& \sum_{k=0}^{Lq} (\epsilon_2 qnt_n^2)^{k} (1-\epsilon_1 qt_n^2)^{n-k} \label{temp_sp}\\
&\leq& \frac{1}{(1-\epsilon_2 qt_n^2)^{Lq}}\sum_{k=0}^{Lq} (\epsilon_2 qnt_n^2)^{k} (1-\epsilon_1 qt_n^2)^{n}  \nonumber\\
&\leq& \frac{1}{(1-\epsilon_2 qt_n^2)^{Lq}}\sum_{k=0}^{Lq} (\epsilon_2 qnt_n^2)^{k} e^{-\epsilon_1 qnt_n^2},\label{sparse_est1}
\end{eqnarray}
where~(\ref{temp_sp}) follows from the fact that~\({n \choose k} \leq n^{k}\) and the estimate~(\ref{sparse_est1}) is true since~\(1-x < e^{-x}\) for all~\(x > 0.\) 

From~(\ref{tn_def}), we have that~\(t_n \rightarrow 0\) as~\(n \rightarrow \infty\) and so using the fact that~\(q\) is a constant, we have that~\((1-\epsilon_2 qt_n^2)^{-Lq} \leq 2\) for all~\(n\) large. Plugging this into~(\ref{sparse_est1}) and using the fact that~\(nt_n^2 \leq \frac{2}{\epsilon_1}\log{n}\) for all~\(n\) large (see~(\ref{tn_def})), we get that
\begin{equation}\label{sparse_est2}
\mathbb{P}\left(\bigcap_{i=1}^{q}E(j_i)\right) \leq 2\sum_{k=0}^{Lq} (\epsilon_2 qnt_n^2)^{k} e^{-\epsilon_1 qnt_n^2} \leq D_1 (\log{n})^{Lq} e^{-\epsilon_1 qnt_n^2},
\end{equation}
for some constant~\(D_1 > 0.\) Setting~\(\alpha_1 = \frac{\alpha}{8}-L\) and again using the expression for~\(t_n\) in~(\ref{tn_def}) we get that the final term in~(\ref{sparse_est2}) is
\begin{equation}
\frac{D_1}{n^{q/8} \cdot (\log{n})^{q\alpha_1}} \cdot \exp\left(-\frac{q\omega_n}{8} + \frac{q\gamma_n}{8}\right) \leq \frac{D_2}{n^{q/8}\cdot (\log{n})^{q\alpha_1}}\exp\left(-\frac{q\omega_n}{8}\right) \nonumber
\end{equation}
for some constant~\(D_2 >0,\)  since~\(\gamma_n <1.\)~\(\qed\)

\subsection*{Left-Right Crossings}
Two~\(t_n \times t_n\) squares~\(S_{i}\) and~\(S_{j}\) are said to be \emph{adjacent} if they share a corner and \emph{plus adjacent} if they share a common side. A sequence of distinct squares~\({\cal Y} = (Y_1,\ldots,Y_w) \subset \{S_j\}\) is said to form a \emph{\(S-\)path} if~\(Y_i\) is adjacent to~\(Y_{i+1}\) for every~\(1 \leq i \leq w-1.\) If, in addition,~\(Y_1\) is also adjacent to~\(Y_w,\) then~\({\cal Y}\) is said to form a~\(S-\)connected cycle. If all the squares in a~\(S-\)path~\({\cal Y}\) are dense, we say that~\({\cal Y}\) is a \emph{dense}~\(S-\)path. Analogous definitions as above hold for~\(S-\)cycles and the plus connected case as well.

\begin{figure}
\centering
\begin{subfigure}{0.5\textwidth}
\centering
\includegraphics[width=3in, trim= 20 150 50 280, clip=true]{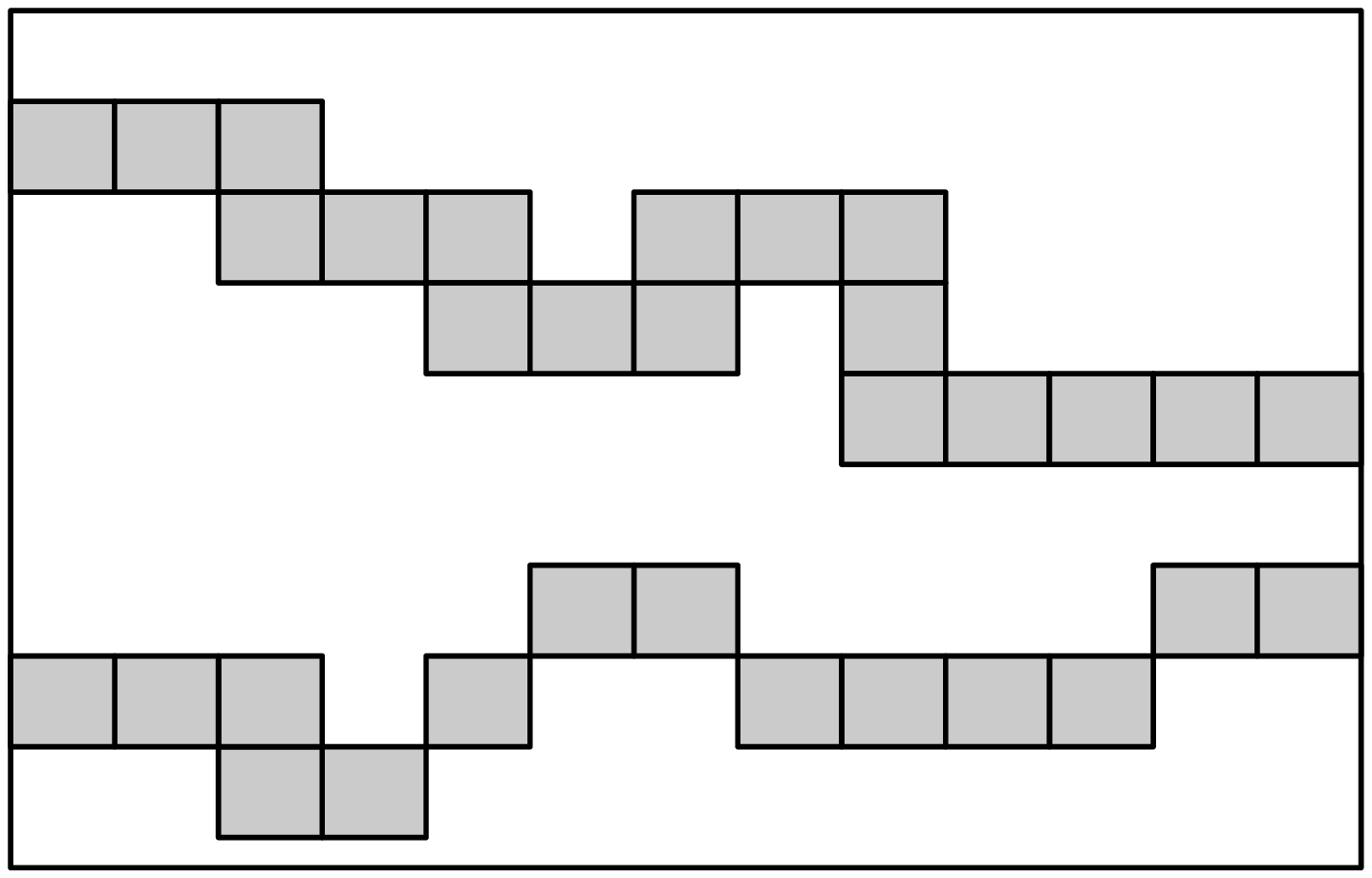}
   \caption{}
\end{subfigure}
\begin{subfigure}{0.5\textwidth}
\centering
   \includegraphics[width=3in, trim= 20 150 50 280, clip=true]{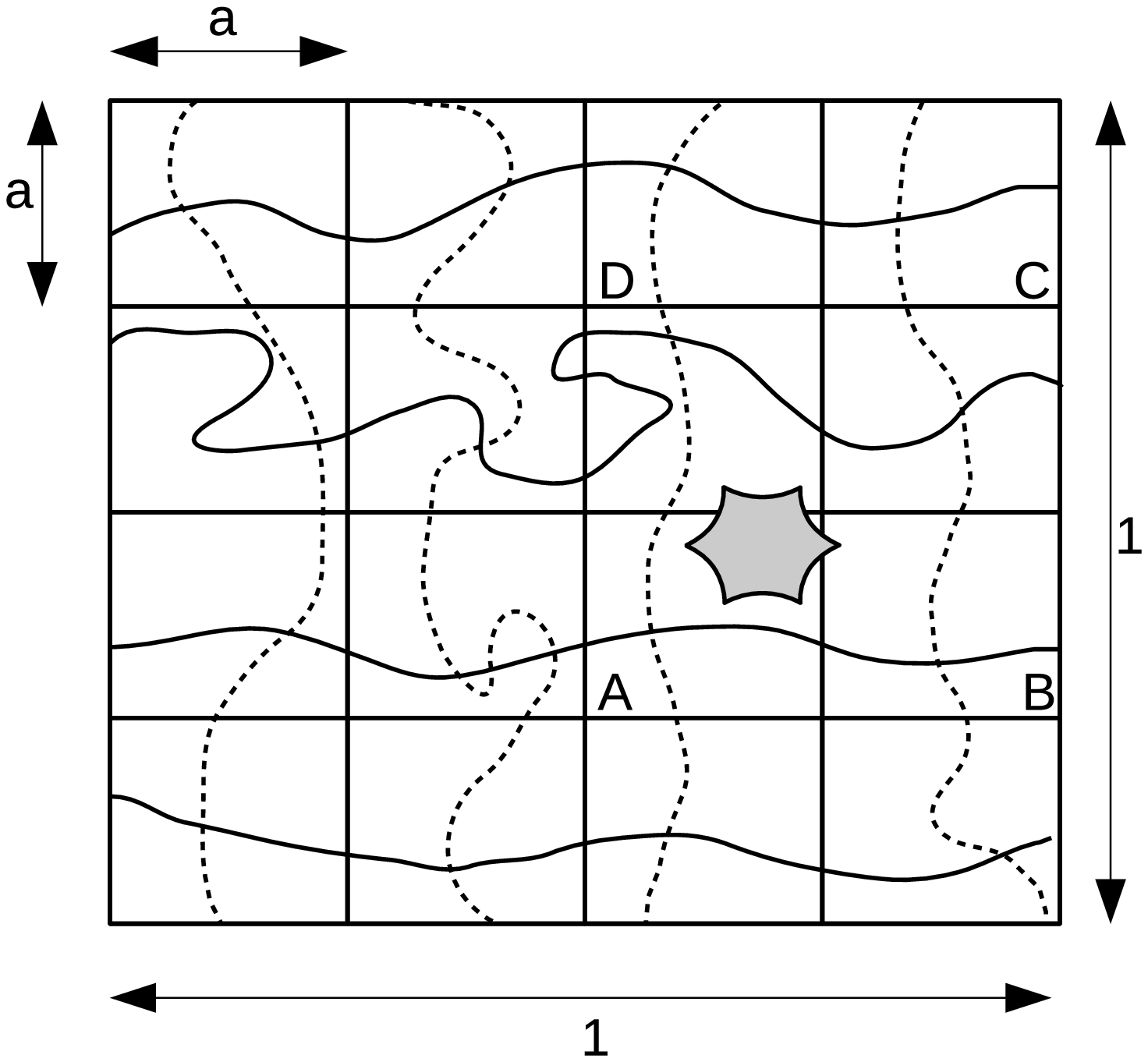}
   \caption{} 
  \end{subfigure}

\caption{\((a)\) Illustration of two left right crossings, where the top left right crossing is plus connected. \((b)\) The backbone formed by the solid and dotted wavy lines due to the occurrence of the event~\(F_n.\) Here~\(a = Mt_n\) is the width of the rectangles and each wavy line is a crossing of the form as shown in~\((a).\) The grey star represents a dense~\(S-\)component distinct from the~\(S-\)component containing the backbone.}
\label{fig_back}
\end{figure}

For future use, we are interested in obtaining a network of long dense~\(S-\)paths that criss-cross each other. We therefore have a couple of additional definitions. For a constant~\(M > 0,\) we divide the unit square~\(S\) into a set of horizontal rectangles \({\cal R}_H\) each of size $1 \times Mt_n$ and also vertically into a set of rectangles \({\cal R}_V,\) each of size $M t_n \times 1.$  If~\((Mt_n)^{-1}\) is an integer, we obtain a perfect tiling as in Figure~\ref{fig_back}\((b).\) Otherwise we choose~\(M_n \in \left[M, M + \sqrt{t_n}\right]\) so that~\((M_nt_n)^{-1}\) is an integer for all~\(n\) large and the tiling of~\(S\) into rectangles of size~$1 \times M_nt_n$ is then perfect. This is possible since~\(t_n \longrightarrow 0\) as~\(n \rightarrow \infty\) (see~(\ref{tn_def})) and so
\[\frac{1}{Mt_n} - \frac{1}{\left(M+\sqrt{t_n}\right)t_n} = \frac{1}{M\sqrt{t_n}} \cdot \frac{1}{M(M+\sqrt{t_n})}\] is bounded below by~\(\frac{1}{2M^2\sqrt{t_n}} \longrightarrow \infty\) as~\(n \rightarrow \infty\)
For notational simplicity, we assume henceforth that~\(M_n = M\) is a constant such that the tiling is perfect.

Let~\(R \in {\cal R}_H \cup {\cal R}_V\) be any rectangle. A distinct sequence of squares\\\({\cal Y} = (Y_1,\ldots,Y_D) \subset \{S_j\}\) contained in~\(R\) is said to be a \emph{left right (top bottom) crossing} of~\(R\) if~\({\cal Y}\) is a~\(S-\)path, the square~\(Y_1\) intersects the left side (top side) of~\(R\) and the square~\(Y_D\) intersects the right side (bottom side) of~\(R.\)  The crossing~\({\cal Y}\) is said to be \emph{dense} if every square in~\(L\) is dense. An analogous definition holds for the plus connected case and for an illustration of left right crossings, we refer to Figure~\ref{fig_back}\((a).\)



For~\(R \in {\cal R}_H,\) let~\(F_n(R)\) be the event that the horizontally long rectangle~\(R \in {\cal R}_H\) contains a dense left right crossing of~\(t_n \times t_n\) squares belonging to \(\{S_j\}_{j \geq 1}.\) Analogously, for \(R \in {\cal R}_V,\) let~\(F_n(R)\) be the event that~\(R\) contains a dense top bottom crossing. Setting
\begin{equation}\label{q_n_def}
F_{n} := \bigcap_{R \in {\cal R}_H \cup {\cal R}_V} F_n(R)
\end{equation}
we have the following result.
\begin{Lemma}\label{fn_lem} We have that
\begin{equation}\label{fn_est}
\mathbb{P}(F_n) \geq 1- \frac{Dn}{n^{M/8}}
\end{equation}
for some constant~\(D>0\) and all~\(n\) large.
\end{Lemma}

\emph{Proof of Lemma~\ref{fn_lem}}: We first prove that
\begin{equation}\label{fn_r_est}
\min_{R \in {\cal R}_{H} \cup {\cal R}_V} \mathbb{P}(F_n(R)) \geq 1 - \frac{D\sqrt{n}}{n^{M/8}}
\end{equation}
for some constant~\(D > 0\) and all~\(n\) large.

To prove~(\ref{fn_r_est}), let~\(R \in {\cal R}_H\) and suppose that~\(F_n(R)\) does not occur.
Necessarily, there exists a sparse plus connected top bottom crossing of~\(R\) (see for example, Theorem~\(2\) of Ganesan~(2017)) and we estimate the probability of such an event happening.  Let~\(\pi = (Y_1,\ldots,Y_q)\) be any plus connected top bottom crossing of~\(R\) and let~\(A_{\pi}\) be the event that every square of~\(\pi\) is sparse. By estimate~(\ref{sj_sparse}) of Lemma~\ref{good_lem}, we have that
\begin{equation}\label{prob_api}
\mathbb{P}(A_{\pi}) \leq \frac{D}{n^{q/8}}
\end{equation}
for some constant~\(D > 0.\)

The square~\(Y_1\) intersects the top edge of~\(R\) and since
\begin{equation}\label{tn_bd}
t_n \geq \frac{r_n}{5} \geq C \cdot \sqrt{\frac{\log{n}}{n}}
\end{equation}
for some constant~\(C >0,\) (see~(\ref{tn_def}) and~(\ref{rn_hom})), there are~\(\frac{1}{t_n} \leq C \sqrt{n}\) choices of~\(Y_1.\) For each such choice of~\(Y_1,\) the number of choices for~\(\pi\) is at most~\(8^{q}\) and since the height of the rectangle~\(R\) is~\(M t_n,\) it is also necessary that~\(q \geq M.\)  Therefore using the geometric summation formula we have that
\begin{equation}\label{fn_te}
\mathbb{P}(F_n^c(R)) \leq  C\sqrt{n} \sum_{q \geq M} 8^{q}\cdot \frac{D}{n^{q/8}}   = CD\sqrt{n} \sum_{q \geq M} \left(\frac{8}{n^{1/8}}\right)^{q}  \leq \frac{D_1\sqrt{n}}{n^{M/8}}
\end{equation}
for some constant~\(D_1 > 0\) and all~\(n\) large. This proves~(\ref{fn_r_est}). 

From~(\ref{fn_r_est}),~(\ref{q_n_def}),~(\ref{fn_te}) and the union bound, we then get
\begin{equation} \label{q_n_est}
\mathbb{P}(F_{n}) \geq 1 - \#({\cal R}_H \cup {\cal R}_V)\frac{D_1\sqrt{n}}{n^{M/8}} \geq 1 - \frac{D_2 n}{n^{M/8}}
\end{equation}
for all~\(n\) large and some constant~\(D_2 > 0,\) where the last estimate in~(\ref{q_n_est}) follows from
\begin{equation}\label{num_rect_est}
\#({\cal R}_H \cup {\cal R}_V) = \frac{2}{Mt_n} \leq D_3\sqrt{n}
\end{equation}
for some constant~\(D_3 > 0,\) by~(\ref{tn_bd}).~\(\qed\)

\setcounter{equation}{0}
\renewcommand\theequation{\thesection.\arabic{equation}}
\section{Proof of Theorem~\ref{thm_rob_rgg}}\label{sec_proof}


\emph{Proof Outline}: Recalling the definition of the event~\(F_n\) in Lemma~\ref{fn_lem}, we set the constant~\( M > 0\) to be sufficiently large so that
\begin{equation}\label{fn_fin_est}
\mathbb{P}(F_n) \geq 1-\frac{1}{n^{10}}.
\end{equation}
If~\(F_{n}\) occurs, then by considering lowermost dense left right crossings of rectangles in~\({\cal R}_H\) and leftmost dense top bottom crossings of rectangles in~\({\cal R}_V,\) we obtain a unique ``backbone" of crossings which we denote as~\({\cal B}.\) This is illustrated in Figure~\ref{fig_back}\((b),\) where the solid and dotted wavy lines together form the backbone~\({\cal B}.\)

Say that a set of squares~\({\cal C} := \{Z_i\}_{1 \leq i \leq r} \subset \{S_j\}\) is a dense~\emph{\(S-\)component} if:\\
\((a)\) For any~\(1 \leq i_1 \neq i_2 \leq r,\) the squares~\(Z_{i_1}\) and~\(Z_{i_2}\) are connected by a dense~\(S-\)path.\\
\((b)\) If~\(Z\) is any dense square adjacent to some square in~\({\cal C},\) then~\(Z \in {\cal C}\) itself.\\
In other words, a dense~\(S-\)component is a \emph{maximal} connected set of dense squares. For a square~\(A \in \{S_j\},\) we define the dense~\(S-\)component~\({\cal T}(A)\) containing~\(A,\) as follows. If~\(A\) is sparse, then define~\({\cal T}(A) = \emptyset;\) else define~\({\cal T}(A)\) to be the dense~\(S-\)component containing square~\(A.\) 

Following the above notation, we let~\({\cal T}({\cal B})\) be the dense~\(S-\)component containing all the squares of~\({\cal B}.\) Using~\({\cal T}({\cal B}),\) we now proceed in three steps to obtain the bridged Hamiltonian cycle. In the first step, we show that with high probability the backbone component~\({\cal T}({\cal B})\) is the only dense component among the squares~\(\{S_j\}_{j \geq 1}.\) In the second step, we show that with high probability, every sparse square is adjacent to some dense square of the backbone component. Finally, we construct the bridged Hamiltonian cycle iteratively by connecting cycles within dense and sparse squares, and then estimate its bridge fraction.


\subsection*{Isolated dense components}
For a square~\(A \in  \{S_j\},\) let
\begin{equation}\label{za_def}
I(A) := F_n \bigcap \{{\cal T}(A) \neq {\cal T}({\cal B})\}
\end{equation}
be the event that the dense component containing the square~\(A\) is not the backbone component~\({\cal T}({\cal B}).\) The event~\(F_n\) guarantees the existence of the backbone and therefore~\({\cal B}\) is well-defined and~\({\cal T}({\cal B}) \neq \emptyset.\) Defining
\begin{equation}\label{ei_def}
I_n := \bigcup_{A \in \{S_j\}} I(A),
\end{equation}
we have that
\begin{equation}\label{ei_main_est}
\mathbb{P}(I_n) \leq \frac{C_I}{(\log{n})^{\alpha-8L+1}}e^{-\frac{3\omega_n}{8}}
\end{equation}
for some constant~\(C_I > 0\) and for all~\(n\) large, where~\(\omega_n \longrightarrow \infty\) is as in~(\ref{rn_hom}).



\emph{Proof of~(\ref{ei_main_est})}: For the~\(t_n \times t_n\) square~\(A \in \{S_j\}\) let~\({\cal U}_1(A)\) be the set of all squares of~\(\{S_j\}\) adjacent to~\(A\) and for~\(i \geq 2,\) let~\({\cal U}_i(A)\) be the set of all squares of~\(\{S_j\}\) adjacent to some square of~\({\cal U}_{i-1}(A)\) so that~\({\cal U}_i(A)\) has~\((2i+1)^2\) squares of~\(\{S_j\}.\)

\begin{figure}[tbp]
\centering
\includegraphics[width=2.8in, trim= 30 290 100 100, clip=true]{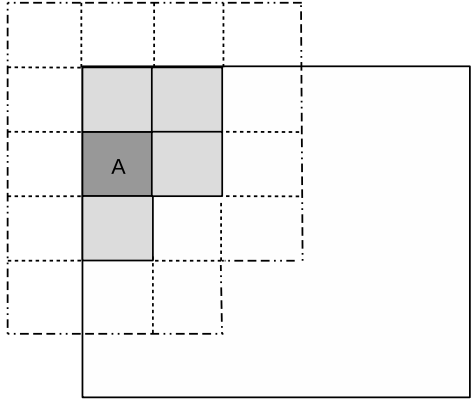}
\caption{The plus connected~\(S-\)cycle~\(L_{cyc}(A)\) shown by the sequence of dotted squares for the square~\(A\) denoted by the dark grey square. The light grey squares together with the square~\(A\) form the dense~\(S-\)component~\({\cal C}(A).\)}
\label{fig_lxly}
\end{figure}

Let~\(S(1+2t_n)\) be the square with same centre as the unit square~\(S\) and of side length~\(1+2t_n.\) Divide the annulus~\(S(1+2t_n)\setminus S\) (of width~\(t_n\)) into~\(t_n \times t_n\) squares~\(\{Q_j\}.\) There necessarily exists a plus connected~\(S-\)cycle~\(L_{cyc}(A) = (Z_1,\ldots, Z_q)\) containing~\(q\) distinct squares in~\(\{S_j\} \cup \{Q_j\}\) surrounding~\({\cal} (A)\) and satisfying the property that~\(Z_1\) is plus adjacent to~\(Z_t\) and~\(Z_i\) is plus adjacent to~\(Z_{i+1}\) for~\(1 \leq i \leq q-1\) (see for example Theorem~\(1\) of Ganesan~(2017)). This is illustrated in Figure~\ref{fig_lxly}, where the dark grey square is~\(A\) and the light and the dark grey squares together form~\({\cal C}(A).\) The sequence of dotted squares are sparse and form~\(L_{cyc}(A).\)

If the event~\(I(A)\) occurs, then~\({\cal C}(A)\) is distinct  from the backbone component and consequently, the plus connected~\(S-\)cycle~\(L_{cyc}(A)\) must be contained in~\({\cal U}_{2M+1}(A).\) Moreover, any square of~\(L_{cyc}(A)\) contained in the interior of the unit square~\(S\) shares a corner with some dense square in~\({\cal C}(A)\) and so is sparse. Thus
\begin{equation}\label{gen_one}
\mathbb{P}(I(A)) \leq \sum_{\pi} \mathbb{P}(L_{cyc}(A)  = \pi),
\end{equation}
where the summation is over all plus connected~\(S-\)cycles contained in~\({\cal U}_{2M+1}(A).\) To evaluate~\(\mathbb{P}(L_{cyc}(A)  = \pi)\) we consider three cases below depending on where the square~\(A\) is located.

\underline{\emph{Case I}}: The square~\(A\) is within a distance~\(2t_n\) from one of the corners of the unit square~\(S.\) There are at least three sparse squares~\(Y_{i_1},Y_{i_2}\) and~\(Y_{i_3}\) of~\(\pi\) that lie in the interior of the unit square~\(S,\) each of which is sparse.  Letting~\(E(i_j)\) denote the event that~\(Y_{i_j}\) is sparse, we get from~(\ref{sj_sparse}) of Lemma~\ref{good_lem} that
\[\mathbb{P}(L_{cyc}(A)  = \pi) \leq \mathbb{P}\left(\bigcap_{j=1}^{3}E(i_j)\right) \leq \frac{C}{n^{3/8} \cdot (\log{n})^{3\alpha_1}}\exp\left(-\frac{3\omega_n}{8}\right)\]
for some constant~\(C > 0,\) where~\(\alpha_1 = \frac{\alpha}{8}-L.\) Plugging this estimate into~(\ref{gen_one}) and using the fact that the number of possibilities for~\(\pi\) depends only on~\(M,\) we get that
\begin{equation}\label{gen_case_one}
\mathbb{P}(I(A))  \leq \frac{D}{n^{3/8} \cdot  (\log{n})^{3\alpha_1}}\exp\left(-\frac{3\omega_n}{8}\right)
\end{equation}
for some constant~\(D > 0.\)

\underline{\emph{Case II}}: The square~\(A\) does not belong to case~\((I)\) but is within a distance of~\(3t_n\) from the boundary of~\(S.\)
In this case at least~\(5\) squares in the \(S-\)cycle~\(L_{cyc}(A)\) lie in the interior of the unit square~\(S.\) Arguing as in Case~\((I)\) above and using~(\ref{sj_sparse}) with~\(q = 5\) gives
\begin{equation}\label{gen_case_two}
\mathbb{P}\left(I(A)\right) \leq  \frac{D}{n^{5/8} \cdot (\log{n})^{5\alpha_1}}\exp\left(-\frac{5\omega_n}{8}\right).
\end{equation}

\underline{\emph{Case III}}: The square~\(A\) is at a distance of~\(3t_n\) away from the boundary of~\(S.\) In this case at least~\(8\) squares in the \(S-\)cycle~\(L_{cyc}(A)\) lie in the interior of the unit square~\(S.\) Arguing as in Case~\((I)\) above and using~(\ref{sj_sparse}) with~\(q = 8\) gives
\begin{equation}\label{gen_case_three}
\mathbb{P}\left(I(A)\right) \leq  \frac{D}{n \cdot (\log{n})^{8\alpha_1}}e^{-\omega_n},
\end{equation}

Next if~\(N_j\) is the number of squares satisfying Case~\((j)\) for~\(j\in \{I, II,III\},\) then we have that
\begin{equation}\label{nj_est}
N_I \leq 16, N_{II} \leq C \cdot \sqrt{n} \text{ and }N_{III} \leq \frac{Cn}{\log{n}}
\end{equation}
for some constant~\(C > 0.\) Indeed, the first estimate on~\(N_I\) is true since there are four corners of~\(S.\) To obtain~\(N_{II},\) use the fact that the number of squares intersecting the boundary of~\(S\) and contained in the interior of~\(S\) is at most~\(\frac{4}{t_n}.\) Therefore the number of squares at a distance of at most~\(3t_n\) from the boundary of~\(S\) is at most~\(\frac{12}{t_n} \leq C\sqrt{n},\) by~(\ref{tn_def}). The final estimate on~\(N_{III}\) is true since the total number of squares in~\(\{S_j\}\) contained in the interior of~\(S\) is~\(\frac{1}{t_n^2} \leq \frac{Cn}{\log{n}},\) again by~(\ref{tn_def}). This proves~(\ref{nj_est}).

Plugging the estimates~(\ref{gen_case_one}),~(\ref{gen_case_two}),~(\ref{gen_case_three}) and~(\ref{nj_est}) into~(\ref{gen_one}) we get
\begin{eqnarray}
\mathbb{P}\left(I_n\right) &\leq& \frac{16 D}{n^{3/8} \cdot (\log{n})^{3\alpha_1}}e^{-\frac{3\omega_n}{8}} + \frac{CD\sqrt{n}}{n^{5/8} \cdot (\log{n})^{5\alpha_1}}e^{-\frac{5\omega_n}{8}} \nonumber\\
&&\;\;\;\;\;\;\;\;\;\;+\;\;\;\frac{Cn}{\log{n}} \cdot \frac{D}{n \cdot (\log{n})^{8\alpha_1}}e^{-\omega_n}. \nonumber
\end{eqnarray}
For all~\(n\) large both~\(n^{3/8}\) and~\(n^{5/8}\) are much larger than any fixed power of~\(\log{n}\) and so we get that
\[\mathbb{P}\left(I_n\right) \leq  \frac{D_1}{(\log{n})^{8\alpha_1+1}}e^{-\frac{3\omega_n}{8}}\]
for all~\(n\) large and some constant~\(D_1 > 0.\) Using~\(\alpha_1 = \frac{\alpha}{8}-L\) we get~(\ref{ei_main_est}).~\(\qed\)


\subsection*{Isolated sparse squares}
Let~\(A \in \{S_j\}\) be any square and let~\(J(A)\) be the event that all the squares adjacent to~\(A\) and contained in the unit square~\(S\) are sparse. Defining
\begin{equation}\label{fi_def}
J_n = \bigcup_{A \in \{S_j\}} J(A)
\end{equation}
we have that
\begin{equation}\label{fi_main_est}
\mathbb{P}(J_n) \leq \frac{C_J}{(\log{n})^{\alpha-8L+1}}e^{-\frac{3\omega_n}{8}}
\end{equation}
for some constant~\(C_J > 0\) and for all~\(n\) large. In particular if the event~\(J^c_n\) occurs, then every sparse square is adjacent to some dense square.\\\\
\emph{Proof of~(\ref{fi_main_est})}: We consider cases~\(I,II\) and~\(III\) as in the previous subsection. In case~\((I)\) there are at least three squares adjacent to~\(A\) and contained in the unit square. Using~(\ref{sj_sparse}) with~\(q = 3\) gives
\begin{equation}\nonumber
\mathbb{P}(J(A)) \leq \frac{C}{n^{3/8} \cdot (\log{n})^{3\alpha_1}}e^{-\frac{3\omega_n}{8}}\text{ for Case I}
\end{equation}
where~\(C > 0\) is a constant. Similarly for case~\((II),\) there are at least~\(5\) squares adjacent to~\(A\) and again using~(\ref{sj_sparse}) with~\(q = 5\) gives
\begin{equation}\nonumber
\mathbb{P}(J(A)) \leq \frac{C}{n^{5/8} \cdot (\log{n})^{5\alpha_1}}e^{-\frac{5\omega_n}{8}} \text{ for Case II.}
\end{equation}
Finally, for case~\((III),\) there are~\(8\) squares adjacent to~\(A\) and so using~(\ref{sj_sparse}) with~\(q = 8\) gives
\begin{equation}\nonumber
\mathbb{P}(J(A)) \leq \frac{C}{n \cdot (\log{n})^{8\alpha_1}}e^{-\omega_n} \text{ for Case III.}
\end{equation}
As before, we let~\(N_j\) be the number of squares in~\(\{S_k\}\) satisfying Case~\((j)\) for~\(j \in \{I,II,III\}.\) Using the estimates for~\(N_I,N_{II}\) and~\(N_{III}\) in~(\ref{nj_est}) and arguing as before, we get~(\ref{fi_main_est}).~\(\qed\)

\subsection*{Constructing the Hamiltonian cycle}
Define the event
\begin{equation}\label{qn_def}
H_n = F_n \bigcap \left(I_n \bigcup J_n\right)^c
\end{equation}
where~\(F_n\) is the ``backbone" event defined in~(\ref{q_n_def}) and the events~\(I_n\) and~\(J_n\) are as in~(\ref{ei_def}) and~(\ref{fi_def}), respectively. From~(\ref{q_n_est}),~(\ref{ei_main_est}) and~(\ref{fi_main_est}), we have that
\begin{equation}\label{hn_est}
\mathbb{P}(H_n) \geq 1-\frac{1}{n^{9}} - \frac{(C_I+C_J)e^{-\frac{3\omega_n}{8}}}{(\log{n})^{\alpha-8L+1}}
\end{equation}
for all~\(n\) large. If the event~\(H_n\) occurs, then there is a backbone~\({\cal B}\) containing dense squares. Recall that~\({\cal T}({\cal B})\) is the dense~\(S-\)component containing all the squares of~\({\cal B}.\) Since the event~\(I^c_n\) occurs, there is no dense star connected~\(S-\)component other than~\({\cal T}({\cal B}).\) Moreover the event~\(J^c_n\) also occurs and so every sparse square is adjacent (i.e. shares a corner) with some dense square in~\({\cal T}({\cal B}).\)

\begin{figure}[tbp]
\centering
\includegraphics[width=2.5in, trim= 130 290 200 230, clip=true]{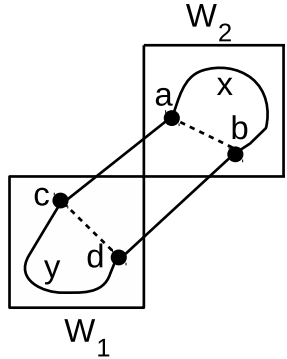}
\caption{Merging the small cycle~\(\tau_1 =\eta_1 = cydc\) contained in the square~\(W_1\) and the small cycle~\(\eta_2 = axba\) contained in the square~\(W_2.\)}
\label{fig_long_rgg}
\end{figure}

We obtain the desired Hamiltonian cycle as follows. Let~\({\cal T}({\cal B}) = \{W_i\}_{1 \leq i \leq t}\) be the set of dense squares in the backbone~\({\cal B}\) and for~\(1 \leq i \leq t,\) let~\(\eta_i\) be the small cycle of edges containing all the nodes of~\(\{X_j\}\) present in the square~\(W_i.\) We now obtain a ``long" cycle~\(\tau({\cal B})\) containing all nodes present in the squares of~\({\cal T}({\cal B})\) as follows.
We set~\(\tau_1 = \eta_1\) and get a series of cycles with increasing lengths, using the small cycles~\(\{\eta_i\}_{1 \leq i \leq t}.\)
The final cycle~\(\tau_t\) would then be the desired long cycle in~\(G.\) First we merge the cycles~\(\tau_1\) and~\(\eta_2\) as in Figure~\ref{fig_long_rgg} by removing one edge each from~\(\tau_1\) and~\(\eta_2\) shown by dotted edges and adding the cross edges shown by straight line segments. The resulting cycle is denoted as~\(\tau_2.\)

To continue the iteration, we now argue that for any~\(2 \leq i \leq t-1,\) the intermediate cycle~\(\tau_i\) still contains an edge of the small cycle~\(\eta_{l(i)}\) contained in the square~\(W_{l(i)},\) adjacent to~\(W_{i+1}.\) This would then allow us to perform the merging between~\(\eta_{i+1}\) and~\(\tau_i\) as described above, and get the cycle~\(\tau_{i+1}.\) The dense square~\(W_{l(i)}\) adjacent to~\(W_{i+1}\) contains at least~\(L \geq 8\) nodes of~\(\{X_j\}\) and so the corresponding small cycle~\(\eta_{l(i)}\) containing all the nodes of~\(W_{l(i)}\) has at least~\(L \geq 8\) edges of~\(G.\) There are exactly~\(8\) squares of~\(\{S_j\}\) adjacent to~\(W_{l(i)}\) and so apart from~\(W_{i+1}\) there are at most~\(7\) squares in the intermediate~\(S-\)component~\({\cal B}_i\) that are adjacent to~\(W_{l(i)}.\) This means that at most~\(7\) edges from the small cycle~\(\eta_{l(i)}\) have been removed so far in the iteration process above.

Continuing this way iteratively for~\(t\) iterations, we get the final cycle~\(\tau_t\) that contains all the nodes of the backbone component. We now iteratively expand the cycle~\(\chi_0 := \tau_t\) by considering sparse squares attached to dense squares in the component~\({\cal T}({\cal B}).\) More precisely, let~\(\{Z_1,\ldots,Z_b\} \subset \{S_j\}\) be the set of all sparse squares. For~\(1 \leq j \leq b,\) let~\(\xi_j, 1 \leq j \leq b\) be a path in~\(G\) containing all the nodes of~\(\{X_k\}\) in the square~\(Z_j.\) As before, we call~\(\{\xi_j\}\) as small paths. Starting from~\(\chi_0,\) we iteratively construct a sequence of intermediate cycles~\(\{\chi_i\}_{1 \leq i \leq b}\) using the paths~\(\{\xi_j\}_{1 \leq j \leq b}.\)

Suppose~\(Z_1\) is adjacent to the square~\(W_{k(1)} \in {\cal T}({\cal B}).\) As argued above, there exists at least~\(L-7\) edges of the small cycle~\(\eta_{k(1)}\) still present in~\(\chi_0.\) Removing one such edge and adding cross edges as in Figure~\ref{fig_long_rgg}, we join the path~\(\xi_1\) and~\(\chi_0\) to get the new cycle~\(\chi_1.\) Repeating the above procedure until all sparse squares are exhausted, we get the final desired Hamilton cycle~\({\cal H}.\)

Summarizing, we begin with~\(\sum_{i=1}^{t} \#\eta_i\) edges belonging to the small cycles and after the iterative procedure described above, we obtain a Hamiltonian cycle containing~\(n\) edges. Therefore if~\(A\) and~\(R\) denote the total number of edges added and removed in the above process, respectively, then
\begin{equation}\label{n_tot}
n = \sum_{i=1}^{t} \#\eta_i - R + A \geq  \sum_{i=1}^{t} \#\eta_i - R
\end{equation}
By definition, the edges in the small cycles belong to the graph~\(G(r_n)\) since~\(t_n < \frac{r_n}{\sqrt{2}}\) (see~(\ref{tn_def})). Therefore it suffices to find an upper bound for~\(A,\) the total number of cross edges added. In each iteration, we remove exactly one edge belonging to the small cycle~\(\eta_i\) in some dense square~\(W_i\) and add two cross edges connecting nodes in~\(W_i\) with nodes in a square adjacent to~\(W_i.\) There are at most eight~\(t_n \times t_n\) squares adjacent to~\(W_i\) and therefore~\(R \leq 8t\) and~\(A \leq 16t.\) Consequently we also get from~(\ref{n_tot}) that
\begin{equation}\label{n_tot2}
n \geq \sum_{i=1}^{t} \#\eta_i - 8t \geq (L-8)t,
\end{equation}
since each dense square contains at least~\(L\) nodes and therefore the small cycle~\(\eta_i\) contains at least~\(L\) edges. Thus from~(\ref{n_tot2}) we see that~\(t \leq \frac{n}{L-8}\) and so the total number of cross edges added is~\(A \leq \frac{16n}{L-8}.\) By construction, the number bridges in the Hamiltonian cycle~\({\cal H}\) is no more than the number of cross edges added and so~\({\cal H}\) has a bridge fraction of at most~\(\frac{16}{L-8}.\)~\(\qed\)





\setcounter{equation}{0}
\renewcommand\theequation{A.\arabic{equation}}
\section*{Appendix}
Writing~\(8nt_n^2 = \theta_n - \gamma_n\) there~\(\theta_n = \frac{1}{\epsilon_1}\left(\log{n} + \alpha\log{\log{n}} + \omega_n\right),\)
it suffices to show that
\[\sqrt{\frac{8n}{\theta_n-1}} - \sqrt{\frac{8n}{\theta_n}} \geq 1.\]
Indeed we have that
\begin{eqnarray}
\sqrt{\frac{8n}{\theta_n-1}} - \sqrt{\frac{8n}{\theta_n}} &=& \sqrt{8n} \left(\frac{\sqrt{\theta_n}-\sqrt{\theta_n-1}}{\sqrt{\theta_n(\theta_n-1)}}\right) \nonumber\\
&\geq& C_1 \sqrt{\frac{n}{\theta_n^2}} \left(\sqrt{\theta_n} - \sqrt{\theta_n-1}\right) \label{gen_a}\\
&=& C_1 \sqrt{\frac{n}{\theta_n^2}} \left(\frac{1}{\sqrt{\theta_n} + \sqrt{\theta_n-1}}\right) \nonumber\\
&\geq& C_2 \sqrt{\frac{n}{\theta_n^2}}  \cdot \frac{1}{\sqrt{\theta_n}} \label{gen_b}\\
&=& C_2 \cdot \sqrt{\frac{n}{\theta_n^3}} \label{gen_c}.
\end{eqnarray}
for some constants~\(C_1,C_2 > 0,\) where~(\ref{gen_a}) and~(\ref{gen_b}) follow from the fact that~\(\theta_n \geq \frac{\log{n}}{\epsilon_1}\) and so~\(\theta_n-1 \geq \frac{\theta_n}{2}\) for all~\(n\) large. Using~\(\theta_n \leq D_2 \log{n}\) for some constant~\(D_2 > 0\) we then get that the final expression in~(\ref{gen_c}) is at least~\(1\) for all~\(n\) large.

\subsection*{Acknowledgement}
I thank Professors Rahul Roy and Federico Camia for crucial comments and for my fellowships.

\bibliographystyle{plain}

\end{document}